\documentclass{article}
\newtheorem{thm}{Theorem}[section]
\newtheorem{cor}[thm]{Corollary}
\newtheorem{prop}[thm]{Proposition}

\newtheorem{lem}[thm]{Lemma}

\newtheorem{ex}{Example}[section]

\newcommand{\be}{\begin{equation}}
\newcommand{\ee}{\end{equation}}
\newcommand{\ben}{\begin{enumerate}}
\newcommand{\een}{\end{enumerate}}
\newcommand{\pa}{{\partial}}

\newcommand{\R}{{\rm R}}
\newcommand{\e}{{\epsilon}}

\newcommand{\pxi}{ {\pa \over \pa x^i}}

\font\BBb=msbm10 at 12pt
\newcommand{\Bbb}[1]{\mbox{\BBb #1}}

\newcommand{\qed}{\hspace*{\fill}Q.E.D.}  
\title{\large \bf Projectively Flat Finsler Metrics of Constant Curvature}
\author{Zhongmin Shen}

\date{Preliminary version in July, revised in September, 2001}

\begin{document}
\maketitle

\begin{abstract}
{It is the Hilbert's Fourth Problem to characterize the (not-necessarily-reversible) distance functions on a bounded convex domain  in $\R^n$ such that   straight lines are shortest paths.  Distance functions induced by a Finsler metric are regarded as smooth ones.  Finsler metrics with straight geodesics said to be projective. 
It is known that the flag curvature of any projective  Finsler metric   is a scalar function of tangent vectors (the flag curvature must be a constant if it is Riemannian). In this paper, 
we study the Hilbert Fourth Problem in the smooth case.  We give a formula for $x$-analytic projective Finsler metrics with constant curvature using a power series with coefficients expressed in terms of $F(0, y)$ and $F_{x^k}(0, y)y^k$. We also give
a formula for general projective Finsler metrics with constant curvature 
using some algebraic equations depending on $F(0, y)$ and $F_{x^k}(0, y)y^k$. 
By these formulas, we obtain several interesting projective Finsler metrics of constant curvature which can be used as models in certain problems.}
\end{abstract}

\section{Introduction}

The well-known  Hilbert's Fourth Problem is to characterize 
the (not-necessarily-reversible)  distance functions on an open subset in $\R^n$ such that straight lines are shortest paths \cite{Hi}. It turns out that there are lots of solutions to the problem.
W. Blaschke first discusses  two-dimensional solutions to the problem \cite{Bl}. 
Later on, R.V. Ambartzumian \cite{Am} and R. Alexander \cite{Al}
independently give all two-dimensional solutions to the problem in a very  elegant way. In the book by A.V. Pogorelov \cite{Po},
he discusses smooth solutions to the problem in three-dimensional case. Z.I. Szabo \cite{Sz} investigates several  problems left by Pogorelov and 
constructs continuous solutions to the problem in high dimensions.  See also \cite{Bu}\cite{AlGeSm} on related issue.

Distance functions induced by a Finsler metrics are regarded as  {\it smooth} ones. The Hilbert Fourth Problem in the smooth case is to characterize
Finsler metrics on an open subset in $\R^n$ whose geodesics are straight lines. Such Finsler metrics are called {\it projective Finsler metrics}.
G. Hamel first characterizes projective Finsler metrics by a system of 
PDE's \cite{Ha}. Later on, A. Rapcs\'{a}k  extends Hamel's result to projectively equivalent Finsler metrics \cite{Rap}. 
It is well-known that every projective Finsler metric
is of scalar curvature, namely, the flag curvature ${\bf K}$ is a scalar function of tangent vectors. It is then natural to determine the  structure of those with constant (flag) curvature.
In the early 20th century,  P. Funk classified all projective Finsler metrics with constant curvature on convex domains in $\R^2$ \cite{Fk1}\cite{Fk2}. Later on, he tried to show  the uniqueness  of 
projectively Finsler metrics with ${\bf K}=1$ on $\Bbb S^2$.
With additional conditions, he showed that the standard Riemannian metric is the only such metric \cite{Fk3}. The final solution is given by R. Bryant
who shows that there is exactly a 2-parameter family of projectively flat Finsler metrics on $\Bbb S^2$ with ${\bf K}=1$ and that the only reversible one is the standard Riemannian metric \cite{Br1}\cite{Br2}. 

In this paper, we first determine the local  structure of 
 $x$-analytic projective Finsler metrics $F(x, y)$ of constant curvature using a power series in $x$.

\begin{thm}\label{thm1}
Let $F(x, y)$ be a projective metric of constant curvature ${\bf K}=\lambda$ on an open neighborhood ${\cal U}$ of the origin in 
$\R^n$. Assume that $F(x, y)$ is $x$-analytic at 
$x=0$. Then $F(x, y)$ is given by
\be 
F(x, y) = \sum_{m=0}^{\infty} {1\over m!}
{d^m\over dt^m} \Big [ \Psi_m (y + tx ) \Big ]|_{t=0},\label{Fxy}
\ee
where $\Psi_m(y)$ is a function on $\R^n$ defined by
\[
\Psi_m :=\cases{ {1\over 2(m+1)} \Big \{
\Big (  \varphi+\psi \Big )^{m+1}- \Big (  \varphi-\psi \Big )^{m+1} \Big \} \ & if $\lambda = -1 $ \cr\\
& \cr \\
\psi \varphi^m \; & if $\lambda =0$ \cr
\\
& \cr \\
{1 \over 2(m+1)i}  \Big \{ 
\Big (   \varphi +i \psi \Big )^{m+1} - \Big (  \varphi-i\psi \Big )^{m+1} \Big \} \ & if $\lambda = 1 $ \cr}
\]
where $\psi(y)$ and $\varphi(y)$ are given by $\psi(y)=F(0, y)$ 
and $\varphi(y)={1\over 2}F(0, y)^{-1}F_{x^k}(0, y) y^k $.
\end{thm}

\bigskip
Note that  a Finsler metric $F(x, y)$ expressed in (\ref{Fxy}) is reversible if and only if 
$\psi(y)$ is reversible ($\psi(-y)=\psi(y)$) and $\varphi(y)$ is anti-reversible ($\varphi(-y)=-\varphi(y)$).
Thus there are lots of reversible projective non-Riemannian Finsler metrics of constant curvature. 
Theorem \ref{thm1} gives us a general formula for 
$x$-analytic projective Finsler metrics of constant curvature. Such Finsler metrics are uniquely determined by $
\psi(y)= F(0, y)$ and $\varphi(y)={1\over 2} F(0, y)^{-1} F_{x^k}(0, y)y^k$. 
Conversely, given  an arbitrary Minkowski norm, $\psi (y)$, and an arbitrary positively homogeneous function of degree one, $\varphi(y)$, on $\R^n$, if the  function $F(x, y)$ defined in (\ref{Fxy}) is
convergent, then it is  a projective Finsler metric of constant curvature ${\bf K}=\lambda$. However, it is 
difficult to determine the domain of convergence around the origin for a given pair, $\{\psi(y), \varphi(y) \}$,  on $\R^n$.

We also study the smooth case.  For any given pair $\{ \psi(y), \varphi(y)\}$, we  
construct a projective  Finsler metric $F(x,y)$ with ${\bf K}= -1, 0$ or $1$, satisfying $F(0, y)=\psi(y)$ and $ {1\over 2}F(0, y)^{-1}F_{x^k}(0, y)y^k =\varphi(y)$.  See Theorems \ref{thmK=-1}, \ref{thmK=0} and \ref{thmK=1} below. 
In this sense,  we have completely determined the local structure of any projective Finsler metrics with constant curvature.  Our method is different from Funk's and Bryant's methods.

Below are some interesting  examples from our constructions. Given a Minkowski norm
$\phi(y)$ on $\R^n$. 
The well-known Funk metric $\Theta(x, y)$ of $\phi(y)$  is defined by
\be
 \phi\Big (y + \Theta(x, y) x \Big ) = \Theta(x, y), \ \ \ \ \ \ y\in T_x{\cal U}, \label{generalFunk}
\ee
where ${\cal U}:= \{ y\in \R^n \ | \ \phi(y) < 1 \}$.
$\Theta(x, y)$ is projective with constant curvature ${\bf K}=-1/4$. 
Moreover,  if the Funk metric $\Theta(x, y)$ is  $x$-analytic at $x=0$, then it can be expressed by
\be
 \Theta (x, y) 
=\sum_{m=0}^{\infty} {1\over (m+1)!} {d^m\over dt^m} \Big [ \phi(y+tx)^{m+1} \Big ]_{t=0}.
\ee
See Example \ref{exFunk} below.

\bigskip

The Funk metric $\Theta(x, y)$ turns out to be a very useful  function. Several important  projective Finsler metrics of constant curvature are constructed using $\Theta(x, y)$. For example,
the Hilbert-Klein metric  on ${\cal U}$ is defined by
\be F(x, y):= {1\over 2}\Big \{ \Theta(x, y)
+ \Theta(x,  -y)\Big \}, \ \ \ \ \ y\in T_x {\cal U}.\label{Hilbertmetric}
\ee
 The Hilbert-Klein metric is a reversible projective Finsler metric with constant curvature ${\bf K}=-1$.  Funk first verified this curvature property for the Hilbert-Klein metric \cite{Fk1}.

We will show that the following Finsler metrics are projective with constant curvature 
${\bf K}= -1$ in its domain,
\be
F(x, y) := {1\over 2} \Big \{ \Theta (x, y)
-\delta \Theta(\delta x, y)\Big \}, \ \  \ \ y\in T_x {\cal U},\label{Shen1}
\ee
where $\delta$ is an arbitrary constant.
When  $-{\cal U}={\cal U}$,
$\Theta(-x, y)=\Theta(x, -y)$. Thus  the Finsler metric  $F$ in (\ref{Shen1}) is just the Hilbert-Klein metric in (\ref{Hilbertmetric}). See Corollary \ref{prop4.1} below.

We will show that the following Finsler metrics are projective with  constant curvature ${\bf K}= -1$ in its domain,
\be
F(x, y):= {1\over 2} \Big \{ \Theta (x, y)
+ {\langle a, y\rangle \over 1 + \langle a, x\rangle } \Big \}, \ \ \ \ \ y\in T_x{\cal U}, \label{GFunk}
\ee
where $a\in \R^n$ is an arbitrary  vector. See Corollary \ref{prop4.2} below.

\bigskip

We will  show that
the  following Finsler metric is projective with constant curvature ${\bf K}=0$ in its domain,
\be
 F(x, y):= \Big \{ 1 +\langle a, x \rangle 
+ {\langle a, y \rangle \over \Theta(x, y)} \Big \} \Big \{ \Theta(x, y)
+ \Theta_{x^i}(x, y)x^i \Big \}, \ \ \ \ \ y\in T_x{\cal U},\label{Shen2}
\ee
 where $a\in \R^n$ is an arbitrary  vector.  This generalizes a result in \cite{Sh4}. See Corollary \ref{prop5.1} below.

We will show that the following Finsler metrics are projective with 
${\bf K}=0$ in its domain,
\be
F(x, y) := { \phi \Big ( (1+\langle a, x\rangle) y - \langle a, y\rangle x \Big ) \over \Big ( 1+ \langle a, x\rangle \Big )^2 }, \ \ \ \ \ y\in T_x{\cal U},
\ee
where $a\in \R^n$ is an arbitrary vector. In fact, this metric is locally Minkowskian. See Corollary \ref{cor6.4}.

\bigskip

Assume that the Funk metric $\Theta(x, y)$ is $x$-analytic so that
 $\Theta(x, y)$ can be extended to an analytic function
$\Theta(z, y)$ in $z\in {\cal U}\otimes \Bbb C \subset \Bbb C^n$. We will show that 
for any angle $\alpha$ with $|\alpha | < \pi/2$, the following Finsler metrics are projective with constant curvature ${\bf K}=1$,
\be
 F(x, y):= {1\over 2}\Big \{ e^{-i\alpha} \Theta (i e^{-i\alpha}x, y)+ \overline{e^{-i\alpha} \Theta (i e^{-i\alpha}x, y)} \Big \}, \ \ \ \ y\in T_x\Bbb R^n. \label{Bryant}
\ee
 When ${\cal U}=\Bbb B^n$, $F(x,y)$ is defined on $\R^n$. It  can be pulled back to the upper and lower semispheres to form
locally projectively flat Finsler metrics on ${\rm S}^n$ with ${\bf K}=1$.
This family of metrics in dimension two are just the Finsler metrics on ${\rm S}^2$ constructed by R. Bryant
\cite{Br1}\cite{Br2}.

\bigskip

The first family of non-projectively flat Finsler metrics of constant curvature
are constructed in \cite{BaSh}. Later on, 
the author constructs infinitely many non-projectively flat Finsler metrics with constant curvature \cite{Sh2}\cite{Sh3}. 
So far, many known non-projectively flat Finsler metrics of constant curvature  are in the form
$F =\alpha+\beta$, where $\alpha$ is a Riemannian metric
and $\beta$ is a $1$-form. Such metrics are called  {\it Randers metrics} \cite{Ran}.  Recently, D. Bao and C. Robles \cite{BaRo} 
have found an equivalent condition for Randers metrics be of constant curvature.

\bigskip
\noindent{\bf 
Note}: After the preliminary version of this paper was sent out,  the author  received an interesting  paper from B. Bryant \cite{Br3}, in which 
Bryant characterizes the (generalized) Finsler metrics on ${\rm S}^n$ with  ${\bf K}=1$ and great circles being as geodesics (see Theorem 2 in \cite{Br3}).   
  As argued by Bryant,  the idea of Theorem 2 in \cite{Br3} can be used to construct projective Finsler metrics with constant curvature ${\bf K}=1$  and a prescribed tangent indicatrix at one point (see Proposition 4 in \cite{Br3}). He also briefly  explains how to construct all of the possible local projective Finsler metrics with ${\bf K}=1$. Bryant informed the author that his idea can also be used to characterize projective Finsler metrics of constant curvature ${\bf K}=0$ or $-1$, although he did not give any detailed discussion on this issue in \cite{Br3}. Bryant's idea \cite{Br1}-\cite{Br3} as well as Funk's idea \cite{Fk1}-\cite{Fk3} are different from ours.  The formula  in (\ref{Fxy}) is not given in Bryant's and Funk's papers.

\bigskip
\noindent{\bf Acknowledgments}: This work was done during author's visit to Prof. S.S. Chern. It was presented at the 2nd national conference on Finsler geometry in Nankai Institute of Mathematics on August 1, 2001.
The author would like to thank S.S. Chern for his warm hospitality and valuable discussion. Thank R. Bryant and Z.I. Szabo for explaining their work to the author.

\section{Preliminaries}

Let $F$  be a  Finsler metric on an $n$-dimensional manifold $M$. 
For a non-zero vector ${\bf y}\in T_pM$, $F$ induces an inner product
$g_{\bf y}$ on $T_pM$ by
\[ g_{\bf y}({\bf u}, {\bf v}) := g_{ij}(x, y) u^iv^j
= {1\over 2} [F^2]_{y^iy^j} (x, y) u^iv^j.\]
Here $x =(x^i)$ denotes the coordinates of $p\in M$ and $(x, y)=(x^i, y^i)$ denotes  the local coordinates of ${\bf y}\in T_pM$.
The geodesics  are characterized by
\[ {d^2c^i\over dt^2}+ 2 G^i \Big (\dot{c}(t) \Big )=0,\]
where $G^i: = {1\over 2} g^{il} \{  [F^2]_{x^ky^l}y^k-[F^2]_{x^l} \}$ are called
the {\it geodesic coefficients} of $F$.
The 
Riemann curvature ${\bf R}_{\bf y}= R^i_{\ k}  dx^k \otimes \pxi|_p : T_pM \to T_pM$ is defined by
\be
R^i_{\ k} = 2 {\pa G^i\over \pa x^k}-y^j{\pa^2 G^i\over \pa x^j\pa y^k}
+2G^j {\pa^2 G^i \over \pa y^j \pa y^k} - {\pa G^i \over \pa y^j}{\pa G^j \over \pa y^k}.  \label{Riemann}
\ee
The Riemann curvature has the following properties: for any non-zero vector ${\bf y}\in T_pM$, 
\[ {\bf R}_{\bf y}({\bf y})=0, \ \ \ \ \ 
g_{\bf y} ({\bf R}_{\bf y}({\bf u}), {\bf v})
=g_{\bf y} ({\bf u}, {\bf R}_{\bf y}({\bf v})),
\ \ \ \ \ {\bf u}, {\bf v}\in T_pM.\]
For a two-dimensional plane $P\subset T_pM$ and a non-zero vector ${\bf y}\in T_pM$, the {\it flag curvature} ${\bf K}(P, y)$ is defined by
\be
{\bf K}(P, {\bf y}):= {g_{\bf y} ({\bf u}, {\bf R}_{\bf y}({\bf u}))
\over g_{\bf y}({\bf y}, {\bf y}) g_{\bf y}({\bf u}, {\bf u})
-g_{\bf y}({\bf y}, {\bf u})^2 },
\ee
where $P={\rm span}\{ {\bf y}, {\bf u}\}$.   
We say that $F$ is of {\it scalar   curvature} ${\bf K}= \lambda(y) $ if 
for any ${\bf y}\in T_pM$, 
the flag curvature ${\bf K}(P, {\bf y})= \lambda({\bf y})$ is independent of $P$ containing ${\bf y}\in T_pM$, that is equivalent to the following system in a local coordinate system $(x^i, y^i)$ in $TM$, 
\[ R^i_{\ k} = \lambda \; F^2 \Big \{ \delta^i_k - F^{-1} F_{y^k} y^i \Big \}.\]
If $\lambda $ is a constant, then $F$ is said to be of {\it constant curvature}.

There are several non-Riemannian quantities in Finsler geometry. 
One of the important non-Riemannian quantities is the E-curvature 
${\bf E}_{\bf y} = E_{ij}  dx^i \otimes dx^j|_p: T_pM \otimes T_pM \to \R$,  defined by
\be
E_{ij}:= {1\over 2} {\pa^2 \over \pa y^i \pa y^j}
\Big [ {\pa G^m \over \pa y^m} \Big ].
\ee
The E-curvature has the following properties: for any non-zero vector ${\bf y}\in T_pM$,
\[ {\bf E}_{\bf y} ({\bf y}, {\bf v}) =0, \ \ \ \ {\bf E}_{\bf y}({\bf u}, {\bf v}) = {\bf E}_{\bf y}({\bf v}, {\bf u}), \ \ \ \ {\bf u}, {\bf v}\in T_pM.\]
For a two-dimensional plane $P\subset T_pM$ and a non-zero vector ${\bf y}\in T_pM$, the {\it flag E-curvature} ${\bf E}(P, y)$ is defined by
\be
{\bf E}(P, {\bf y}):= { F^3({\bf y})  {\bf E}_{\bf y} ({\bf u}, {\bf u})
\over g_{\bf y}(y, y) g_{\bf y}({\bf u}, {\bf u})
-g_{\bf y}({\bf y}, {\bf u})^2 },
\ee
where $P={\rm span}\{ {\bf y}, {\bf u}\}$.  
We say that $F$ has constant E-curvature ${\bf E}= (n+1)c$ if 
for any pair $(P, {\bf y})$, ${\bf E}(P, {\bf y}) = (n+1) c$, that is equivalent to the following system, 
\[ E_{ij} = (n+1) c F_{y^iy^j}.\]
We know that the Funk metric on a strongly convex domain 
satisfies that ${\bf K}=-1/4$ and ${\bf E}= (n+1)/4$ \cite{Sh1}.

\section{Projective Finsler Metrics}

A distance function  on a set ${\cal U}$ is a function 
$d: {\cal U}\times {\cal U}\to \R$ with the following properties 
\ben
\item[(a)] $d(p, q)\geq 0$ and equality holds if and only if $p=q$;
\item[(b)] $d(p, q)\leq d(p, r)+ d(r, q)$.
\een
A distance function on a convex domain ${\cal U}\subset \R^n$ is said to be {\it projective} (or {\it rectilinear}) if 
straight lines are shortest paths. The Hilbert's Fourth Problem is to characterize projective distance functions. 

\bigskip

A  distance function $d$ on a manifold $M$ is said to be {\it smooth}  if it is induced by a Finsler metric $F$ on $M$,
\[ d(p, q):=\inf_{c} \int_0^1 F(\dot{c}(t)) dt,\]
where the infimum is taken over all curves $c(t)$, $0 \leq t \leq 1$,  joining $p=c(0)$ to $q=c(1)$. Thus smooth  distance functions can be studied using calculus (see \cite{AIM}, \cite{BCS} and \cite{Sh1}, etc.).

Now we start to discuss smooth projective distance functions, or projective Finsler metrics on an open domain ${\cal U}\subset \R^n$.
First, let us use the following notations.  
The local coordinates of a tangent vector ${\bf y}= y^i \pxi|_p\in T_x{\cal U}$ will be denoted by $(x, y)$. Hence all quantities are functions of $(x, y)\in {\cal U}\times \R^n$.
It is known that 
a Finsler metric  $F(x,y)$ on  ${\cal U}$  is  projective  if and only if 
its geodesic coefficients $G^i$ are in the form 
\[G^i(x, y) = P(x, y) y^i,\]
 where $P: T{\cal U} = {\cal U}\times \R^n \to \R$ is positively homogeneous with degree one, 
$P(x,  \lambda y) = \lambda P(x, y)$, $\lambda >0$. We call $P(x, y)$ the {\it projective factor} of $F(x, y)$.
The following lemma plays an important role.

\begin{lem}\label{lemRap}{\rm (Rapcs\'{a}k \cite{Rap})} Let $F(x, y)$ be a Finsler metric on an open subset ${\cal U} \subset \R^n$. $F(x, y)$ is projective on ${\cal U}$ if and only if 
it satisfies 
\be
F_{x^k y^l} y^k = F_{x^l}. \label{Rap}
\ee
In this case, the projective factor $P(x, y)$ is given by
\be
P = { F_{x^k} y^k\over 2 F}.\label{Gi}
\ee
\end{lem}

\bigskip
Much earlier, G. Hamel \cite{Ha} proved that a Finsler metric $F(x, y)$ 
on ${\cal U}\subset \R^n$ is projective   if and only if 
\be
F_{x^ky^l} = F_{x^l y^k}.\label{Ham}
\ee
Thus (\ref{Gi}) and (\ref{Rap}) are equivalent.

\bigskip

Let $F(x,y)$ be a projective Finsler metric on ${\cal U}\subset \R^n$ and $P(x,y)$ its projective factor.
Put 
\be
 \Xi:= P^2 - P_{x^k}y^k.\label{Xi}
\ee
 Plugging $G^i=Py^i$ into (\ref{Riemann}) yields
\be
R^i_{\ k} = \Xi  \; \delta^i_k +\tau_k
 \; y^i,\label{Rik}
\ee
where
\be
\tau_k =  3 (P_{x^k} - P P_{y^k} ) + \Xi_{y^k}.\label{tau}
\ee
See \cite{Sh1} for more discussion.
It follows from (\ref{Rik}) and (\ref{tau}) that the Ricci curvature ${\bf Ric}:=R^k_{\ k}$
is given by 
\be
 {\bf Ric} = (n-1) \Xi. \label{Ric}
\ee
By the symmetry property that
$g_{ji} R^i_{\ k} =g_{ki}R^i_{\ j}$, we can show that 
\be
R^i_{\ k}  
= \Xi \Big \{ \delta^i_k - F^{-1}F_{y^k} y^i \Big \}. \label{Ri_kRi_k}
\ee
Comparing
(\ref{tau}) and (\ref{Ri_kRi_k}), we  obtain 
\be
P_{x^k} - PP_{y^k} = - { (\Xi F)_{y^k} \over 3 F }.\label{PXi}
\ee

\bigskip

From (\ref{Ric}) and (\ref{Ri_kRi_k}), we immediately obtain the following
\begin{lem}\label{prop2.2} 
For a locally projectively flat Finsler metric $F$ on an $n$-manifold $M$,
the flag curvature and the Ricci curvature are related by 
\be
{\bf K}(P, {\bf y}) ={1\over n-1} { {\bf Ric}({\bf y}) \over F^2({\bf y}) }, \ \ \ \ \ \ {\bf y}\in P\subset T_pM.
\ee
\end{lem}

It follows from  Lemma \ref{prop2.2} that  a locally projectively flat Finsler metric  has constant Ricci curvature if and only if it has constant flag curvature.

\bigskip

The following is our key lemma to determine 
the local metric structure of projective Finsler metrics with constant curvature.

\begin{lem}\label{lemkey} 
Let $F(x,y)$ be a Finsler metric on an open subset ${\cal U}\subset \R^n$.
Then $F(x,y)$ is projective if and only if  there is a positively homogeneous function with degree one,
$P(x, y)$,   and a positively  homogeneous function of degree zero, $\lambda(x, y)$, on $T{\cal U}={\cal U}\times \R^n$ such that 
\be
F_{x^k} = (P F)_{y^k} \label{eq1}
\ee
\be
P_{x^k} = PP_{y^k} - {1\over 3F} ( \lambda F^3 )_{y^k}.\label{eq2}
\ee
In this case, $P= {1\over 2} F^{-1} F_{x^k}y^k$ and 
$F$ is of scalar curvature ${\bf K}=\lambda$.
\end{lem}
{\it Proof}: Assume that $F$ is projective. 
By Lemma \ref{lemRap}, $F=F(x, y)$ satisfies (\ref{Rap}) and  the geodesic coefficients  are in the form $G^i = P y^i$,  where 
\[ P := {F_{x^l} y^l\over 2F}.\]
$P= P(x, y)$ is a positively homogeneous function of degree one on $T{\cal U}={|cal U}\times \R^n$.
Observe that 
\[
(PF)_{y^k} = {1\over 2} (F_{x^l}y^l)_{y^k}
= {1\over 2} (F_{x^ly^k}y^l + F_{x^k})
= {1\over 2} (F_{x^k} + F_{x^k}) = F_{x^k}.\]
Thus $F$ satisfies (\ref{eq1}). 
Let 
\[ \lambda: ={\Xi \over F^2},\]
where $\Xi := P^2-P_{x^k}y^k$. $\lambda =\lambda(x, y)$ is a positively homogeneous function of degree zero on $T{\cal U}={\cal U}\times \R^n$.
Plugging it  into 
(\ref{PXi}) yields
\[ P_{x^k} - P P_{y^k} = -{1\over 3F} (\lambda  F^3)_{y^k}.\]
Thus $P$ satisfies (\ref{eq2}).

Suppose that there are a positively homogeneous function of degree one,
$P(x,y)$,  and a positively  homogeneous function of degree zero, $\lambda(x, y)$, on $T{\cal U}$ such that  (\ref{eq1}) and (\ref{eq2}) hold.
First by (\ref{eq1}), we obtain
\[
F_{x^ky^l}y^k = (PF)_{y^ky^l} y^k
= (PF)_{y^l} = F_{x^l},\]
and
\[ P = { (PF)_{y^k}y^k\over 2F} = { F_{x^k} y^k\over 2F}.\]
By Lemma \ref{lemRap}, we conclude that $F$  is projective  with  the geodesic coefficients in the form $G^i = P y^i$. Plugging $G^i=Py^i$ into  (\ref{Riemann}) yields a formula for ${\bf Ric} = R^m_{\ m}$
\[ {\bf Ric} = (n-1)\Xi,\]
where $\Xi:= P^2-P_{x^k}y^k$. 
Contracting (\ref{eq2}) with $y^k$ yields  that 
\[{\bf Ric}  =(n-1) ( P^2 - P_{x^k}y^k )= 
 { (n-1) \over 3F} (\lambda F^3)_{y^k} y^k =  (n-1) \lambda F^2.\]
By Lemma \ref{prop2.2}, $F$  has flag curvature ${\bf K} = \lambda$. 
\qed

\section{Solving the inverse problem}

In this section, we are going to determine the local metric structures of $x$-analytic projective metrics $F(x, y)$ with constant curvature by solving the first order  partial differential equations (\ref{eq1}) and (\ref{eq2}).
Here a projective Finsler  metric $F(x, y)$
on ${\cal U}$ is said to be {\it $x$-analytic} at $x_o=(x^i_o)\in {\cal U}$, if
there is a number $\e >0$ such that $F(x, y)$ can be expressed as a power series for $x=(x^i)\in {\cal U}$ with  $|x-x_o| <\e$,
\[ F(x, y) = \sum_{m=0}^{\infty} {1\over m!}\sum_{i_1\cdots i_m=1}^n a_{i_1\cdots i_m}(y) (x^{i_1}-x^{i_1}_o) \cdots (x^{i_m}-x^{i_m}_o),\]
where $a_{i_1\cdots i_m}(y)$ are $C^{\infty}$ functions on $\R^n\setminus\{0\}$ with $a_{i_1\cdots i_m}(\lambda y)
=\lambda a_{i_1\cdots i_m} (y)$, $\lambda >0$.

\bigskip

Let $F$ be a Finsler metric defined on an open neighborhood ${\cal U}$ of the origin 
in $\R^n$. Assume that 
$F(x, y)$ is projective with  constant curvature ${\bf K}=\lambda$ on an open neighborhood of the origin in $\R^n$.
By Lemma \ref{lemkey}, there is a positively homogeneous function of degree one,
$P(x, y)$  on $T{\cal U}= {\cal U}\times \R^n$ such that (\ref{eq1}) and (\ref{eq2}) hold,
\be
F_{x^k} = (PF)_{y^k}, \ \ \ \ \ P_{x^k} = P P_{y^k} -\lambda F F_{y^k}.
\label{FPFP}
\ee

Define $\psi_m(x, y)$ by
\[
\psi_m :=\cases{ {1 \over 2(m+1)\e} \Big \{
\Big ( P+\e F \Big )^{m+1} - \Big ( P-\e F \Big )^{m+1} \Big \} \ & if $\lambda = -\e^2 $ \cr\\
& \cr \\
F P^m \; & if $\lambda =0$ \cr
\\
& \cr \\
{1 \over 2(m+1)\e i}  \Big \{ 
\Big ( P+i\e F \Big )^{m+1} - \Big ( P-i \e F\Big )^{m+1} \Big \} \ & if $\lambda =  \e^2 $ \cr}
\]
where $\e >0$ is an arbitrary constant.
We have
\[ \psi_0 = F, \ \ \ \ \psi_1 = P F, \ \ \ \
\psi_2 = F \Big ( P^2 - {\lambda \over 3} F^2 \Big ) , \ \ \cdots.\]

Let us first consider the case when $\lambda =-\e^2 >0$.
Using (\ref{FPFP}), we obtain
\begin{eqnarray*}
&&\hspace{-1 cm}  \Big [ {1\over m+1}\Big (P+\e F\Big )^{m+1}\Big ]_{x^k}\\
&& = \Big (P+\e F\Big )^m \Big (  P_{x^k}+ \e F_{x^k}   \Big )\\
& & = \Big (  P+\e F \Big )^m 
\Big (  P P_{y^k} + \e^2 F F_{y^k} +\e (PF)_{y^k} \Big )\\
&& = \Big ( P+\e F \Big )^{m+1} 
\Big ( P+\e F \Big )_{y^k}\\
&& = {1 \over m+2} \Big [ \Big ( P+\e F \Big )^{m+2}\Big ]_{y^k}.
\end{eqnarray*}
Similarly,
\[ \Big [{1 \over m+1} \Big (P-\e F \Big )^{m+1}\Big ]_{x^k}
= {1 \over m+2}\Big [ \Big (P-\e F \Big )^{m+2}\Big ]_{y^k}.
\] 
By the above two equations, we obtain 
\be
\Big [\psi_m\Big ]_{x^k} 
=\Big  [\psi_{m+1}\Big ]_{y^k}.
\ee
Thus
\be
F_{x^{i_1}\cdots x^{i_m} } = \Big [ \psi_0\Big ]_{x^{i_1}\cdots x^{i_m}}
= \Big [ \psi_1 \Big ]_{x^{i_2}\cdots x^{i_m}y^{i_1} }
= \cdots = \Big [ \psi_m \Big ]_{y^{i_1}\cdots y^{i_m} }. \label{eq26}
\ee
One can verify that (\ref{eq26}) holds for other cases when $\lambda=0$ and $\e^2$.

Let 
\[ \psi(y):= F(0, y), \ \ \ \ \ \ \varphi(y):= P(0, y).\]
Define $\Psi_m(y)$ by
\be
\Psi_m :=\cases{ {1 \over 2(m+1)\e} \Big \{
\Big (  \varphi+\e \psi  \Big )^{m+1} - \Big (   \varphi -\e \psi \Big )^{m+1} \Big \} & if $\lambda = - \e^2 $ \cr\\
& \cr \\
\psi \varphi^m \; & if $\lambda =0$ \cr
\\
& \cr \\
{1\over 2(m+1)\e i }  \Big \{
\Big ( \varphi + i \e \psi  \Big )^{m+1} - \Big (   \varphi-i\e \psi  \Big )^{m+1} \Big \} & if $\lambda =  \e^2 $ \cr}\label{Psi_m}
\ee
We see that $\psi_m(0, y) = \Psi_m (y)$. Setting $x=0$ in
(\ref{eq26}) yields 
\[ F_{x^{i_1}\cdots x^{i_m} }(0, y) =\Big [\Psi_m\Big ]_{y^{i_1}\cdots y^{i_m}}(y).\]
We obtain
\[\sum_{i_1\cdots i_m=1}^n  F_{x^{i_1}\cdots x^{i_m} }(0, y) x^{i_1} \cdots x^{i_m}
= {d^m\over dt^m} \Big [ \Psi_m (y+t x) \Big ]_{t=0}.
\]
By assumption that $F(x, y)$ is $x$-analytic at $x =0$, we obtain
\be
 F(x, y) = \sum_{m=0}^{\infty} {1\over m!} {d^m\over dt^m}
\Big [ \Psi_m (y+tx) \Big ]_{t=0}.\label{thm2F}
\ee

\bigskip

Define $\Phi_m(y)$ by 
\be
\Phi_m :=\cases{ {1 \over 2(m+1)} \Big \{
\Big (   \varphi + \e \psi \Big )^{m+1} +  \Big (   \varphi  -\e \psi \Big )^{m+1} \Big \} \ & if $\lambda = -\e^2 $ \cr\\
& \cr \\
{1\over m+1} \varphi^{m+1} \; & if $\lambda =0$ \cr
\\
& \cr \\
{1 \over 2(m+1)}  \Big \{
\Big (  \varphi + i \e \psi \Big )^{m+1} + \Big (   \varphi- i \e \psi \Big )^{m+1} \Big \} \ & if $\lambda =  \e^2 $ \cr}\label{Phi_m}
\ee
By a similar argument, we obtain 
\[ P_{x^{i_1}\cdots x^{i_m}} (0, y)
= \Big [ \Phi_m \Big ]_{y^{i_1}\cdots y^{i_m}} (y).\]
Since $F(x, y)$ is $x$-analytic at $x=0$, the
projective factor  $P(x, y) := {1\over 2}F^{-1}F_{x^k} y^k$ is $x$-analytic at $x=0$ too. We obtain \be
P(x, y) 
= \sum_{m=0}^{\infty} {1\over m!} {d^m \over dt^m} \Big [ \Phi_m ( y
+ t x ) \Big ]|_{|t=0}. \label{thm2P}
\ee
We have proved the following
\begin{thm}\label{thm2} 
Let $F(x, y)$ be a projective Finsler metric of constant curvature ${\bf K} 
= \lambda $ on an open subset ${\cal U}\subset \R^n$. Suppose that 
$F(x, y)$ is $x$-analytic at $x=0$. Then 
$F$ and $P:={1\over 2} F^{-1}F_{x^k}y^k$  are given by 
(\ref{thm2F}) and (\ref{thm2P}) respectively.
\end{thm}

\bigskip

\begin{ex}\label{exFunk}
{\rm (Funk metric)
Let
$\phi(y)$ be a Minkowski norm on $\R^n$
and 
\[ {\cal U}:= \Big \{ y\in \R^n \ | \ \phi(y) < 1 \Big \}.\]
Define  $\Theta(x, y) >0$  by
\be
\Theta(x, y)= \phi \Big (y+\Theta(x, y) x \Big ), \ \ \ \ \ y\in T_x{\cal U}.\label{FunkFunkm}
\ee
 $\Theta(x, y)$ is a Finsler metric on ${\cal U}$ satisfying
the following system of equations \cite{Ok},
\be
\Theta_{x^k} = \Theta\Theta_{y^k}. \label{Funk}
\ee
The function $\Theta(x, y)$ is called the {\it Funk metric} of $\phi(y)$ on ${\cal U}$.
Let
\[ F(x, y):= \Theta(x, y), \ \ \ \ \ \ P(x, y):= {1\over 2} \Theta(x, y).\]
 Note that $F(0, y) =\phi(y)$ and $P(0, y)= {1\over 2} \phi(y)$. 
We see that (\ref{Funk}) is equivalent to (\ref{FPFP}) with $\lambda = -1/4$. Thus $F(x, y)$ is a projective Finsler metric with ${\bf K}=-1/4$
and its projective factor $P(x, y)= {1\over 2}F(x, y)$. 

Observe that $\varphi(y) + {1\over 2} \psi(y)= \phi(y)$. Thus 
\[ \Psi_m (y) = {1\over m+1} \phi(y)^{m+1}.\]
Assume that $\Theta(x, y)$ is $x$-analytic  at $x=0$. Then 
$\Theta(x, y)$ can be expressed by 
\be
\Theta(x, y) = \sum_{m=0}^{\infty}
{1\over (m+1)!}  {d^m \over dt^m} \Big [ \phi(y+tx)^{m+1} \Big ]|_{t=0}.\label{Funkmetric}
\ee
}
\end{ex}

\section{${\bf K}=-1$}
In this section, with the above  discussion on analytic solutions, we are going to construct smooth solutions $F(x, y)$ with $F(0, y)=\e\psi(y)$ and $F_{x^k}(0, y)y^k= 2\psi(y) \varphi(y)$ for any given pair $\{ \psi(y), \varphi(y)\}$ and a positive number $\e >0$.

Let $F(x,y)$ be a projective Finsler metric of constant curvature 
${\bf K}=-1$ on an open neighborhood ${\cal U}$ of the origin in $\R^n$. 
Let $\e >0$ be an arbitrary constant. Then $F_{\e}(x, y):=\e^{-1}F(x, y)$ has constant curvature ${\bf K}=-\e^2$.
Suppose that $F(x, y)$ is $x$-analytic at $x=0$. According to Theorem \ref{thm2}, $F_{\e}(x, y)$ and 
$P_{\e}(x, y):={1\over 2} (F_{\e})^{-1} (F_{\e})_{x^k}y^k$ are given by 
\begin{eqnarray}
F_{\e}(x, y) & = &  \sum_{m=0}^{\infty} {1\over m!}
{d^m\over dt^m} \Big [ \Psi_m(y+tx ) \Big ]|_{t=0},\label{FK=-1}\\
P_{\e}(x, y) & = & \sum_{m=0}^{\infty} {1\over m!}
{d^m\over dt^m} \Big [ \Phi_m(y+tx ) \Big ]|_{t=0}, \label{PK=-1}
\end{eqnarray}
where $\Psi_m(y)$ and $\Phi_m(y)$ are given by
\begin{eqnarray*}
\Psi_m  :& = &  {1\over 2 (m+1)\e} \Big \{ 
\Big (  \varphi+\e \psi \Big )^{m+1}
- \Big (  \varphi -\e \psi \Big )^{m+1}\Big \},\\
\Phi_m :& = &  {1\over 2 (m+1)} \Big \{ 
\Big (   \varphi+\e \psi \Big )^{m+1}
+ \Big ( \varphi -\e \psi \Big )^{m+1}\Big \},
\end{eqnarray*}
where $\psi (y)=F_{\e}(0, y)$  and $\varphi(y)=P_{\e}(0, y)$.
\bigskip

For $y\not=0$, let $\Psi_{\e}=\Psi_{\e}(x, y)$ and $\bar{\Psi}_{\e}=\bar{\Psi}_{\e}(x, y)$ be solutions to the following equations
by
\be
 \varphi ( y + \Psi_{\e} \; x )+\e \psi ( y + \Psi_{\e} \; x )  =\Psi_{\e},\label{p1}
\ee
\be
 \varphi ( y + \bar{\Psi}_{\e} \; x )-\e \psi ( y + \bar{\Psi}_{\e} \; x )  =\bar{\Psi}_{\e}.\label{p2}
 \ee 
There is a constant $C>0$ such that for any $x$ with $(1+\e)|x| < C$ and any $y \in \R^n$, the above two systems have unique solutions $\Psi_{\e}(x,y)$ and $\bar{\Psi}_{\e}(x, y)$, respectively. 
By (\ref{p1}) and (\ref{p2}), one can show that the functions  $\Psi_{\e}$ and $\bar{\Psi}_{\e}$ satisfy the following equations,
\be
(\Psi_{\e})_{x^k} = \Psi_{\e}(\Psi_{\e})_{y^k},
\ \ \ \ \ \ (\bar{\Psi}_{\e})_{x^k} = \bar{\Psi}_{\e} (\bar{\Psi}_{\e})_{y^k}.\label{Funkeq1}
\ee
By (\ref{Funkeq1}), we obtain
\be
(\Psi_{\e})_{x^{i_1}\cdots x^{i_m}}(x, y)
= {1\over m+1} \Big [ (\Psi_{\e})^{m+1} \Big ] _{y^{i_1}\cdots y^{i_m} } (x, y). \label{FunkA}
\ee
Note that $\Psi_{\e}(0, y) = \varphi(y) +\e \psi(y)$. Setting $x=0$ in 
(\ref{FunkA}) yields 
\[ (\Psi_{\e})_{x^{i_1}\cdots x^{i_m}}(0, y)
= {1\over m+1} \Big [ \Big ( \varphi+\e \psi \Big )^{m+1} \Big ] _{y^{i_1}\cdots y^{i_m} } (y).\]
We assume that $\Psi_{\e} (x, y)$ is $x$-analytic at $x=0$.
By th above identities,  we obtain
\[ \Psi_{\e}(x, y)= \sum_{m=0}^{\infty} 
{1\over (m+1)!} {d^m\over dt^m}\Big [ \Big ( \varphi(y+ t x)+\e \psi(y+ t x)\Big )^{m+1} \Big ]|_{t=0}.\]
We assume that $\bar{\Psi}_{\e}(x, y)$ is $x$-analytic at $x=0$. 
By a similar argument, we obtain
\[
\bar{\Psi}_{\e}(x, y)=  \sum_{m=0}^{\infty} 
{1\over (m+1)!} {d^m\over dt^m}\Big [ \Big (\varphi(y+ t x)-\e\psi(y+ t x) \Big )^{m+1} \Big ]|_{t=0}.
\]
By the above power series, we can express $F_{\e}(x, y)$ in (\ref{FK=-1}) 
and $P_{\e}(x, y)$ in (\ref{PK=-1}) by
\begin{eqnarray*}
F(x, y)& = &  \e F_{\e}(x, y)= {1\over 2}\Big \{ \Psi_{\e}(x, y)
- \bar{\Psi}_{\e}(x, y) \Big \}\\
P(x, y) & = &  P_{\e}(x, y)= {1 \over 2}\Big \{ \Psi_{\e}( x, y)
+ \bar{\Psi}_{\e}(x, y) \Big \}.
\end{eqnarray*}
with $F(0, y) = \e \psi(y)$ and $ P(0, y) = \varphi(y).$

\bigskip

The above arguments lead to the following

\begin{thm}\label{thmK=-1} Let $\psi=\psi(y)$ be a Minkowski norm on $\R^n$ and $\varphi:=\varphi(y)$ be a positively homogeneous function of degree one on $\R^n$. For an arbitrary constant $\e>0$, let $\Psi_{\e}=\Psi_{\e}(x, y)$ and $\bar{\Psi}_{\e}=\bar{\Psi}_{\e}(x, y)$ denote the functions defined by 
(\ref{p1}) and (\ref{p2}) respectively. Then the following function
\be
F(x, y):=  {1\over 2}\Big \{ \Psi_{\e}(x, y)
- \bar{\Psi}_{\e}(x, y) \Big \} \label{eq32}
\ee
is a projective Finsler metric on its domain with constant curvature ${\bf K}=-1$ and its projective factor $P(x, y)= {1\over 2} F^{-1} F_{x^k}y^k$ is given by 
\be
P(x, y): = {1 \over 2}\Big \{ \Psi_{\e}( x, y)
+ \bar{\Psi}_{\e}(x, y) \Big \}.\label{eq33}
\ee
Further, $F(0, y)=\e \psi(y)$ and $P(0, y)=\varphi(y)$.
\end{thm}
{\it Proof}: It follows from  (\ref{Funkeq1}) that $F$ satisfies (\ref{Rap}) and hence it is projective. Observe that 
\[ F_{x^k} y^k
= {1\over 4} \Big [ \Psi_{\e}^2 + \bar{\Psi}_{\e}^2 \Big ]_{y^k} y^k
= {1\over 2} \Big \{ \Psi_{\e}^2 + \bar{\Psi}_{\e}^2 \Big \}.\]
Thus the projective factor $P= {1\over 2} F^{-1}F_{x^k}y^k$ is given by
\[ P = {1\over 2} { \Psi_{\e}^2 - \bar{\Psi}_{\e}^2 \over  \Psi_{\e}+\bar{\Psi}_{\e} } = {1\over 2} \{  \Psi_{\e}- \bar{\Psi}_{\e} \Big \}.\] 
That is, $P$ satisfies (\ref{eq33}).
By a similar argument, we obtain
\[ \Xi = P^2 - P_{x^k}y^k = - F^2.\]
Thus the flag curvature ${\bf K}=-1$ by Lemma \ref{prop2.2}.
\qed

\bigskip

Let $\phi(y)$ be an arbitrary Minkowski norm on $\R^n$
and ${\cal U}:= \{ y\in \R^n \ | \ \phi(y) < 1 \}.$
Let $\Theta(x, y)$ denote the Funk metric of $\phi$ on ${\cal U}$ defined by (\ref{FunkFunkm}).
Let $\e >0$ and  $\delta$ be  arbitrary positive constants and 
\[ \psi(y) := {1-\delta \over 2\e}  \phi(y), \ \ \ \ \ \varphi(y):= {1+\delta \over 2} \phi(y).\]
We have
\[ \varphi+\e \psi = \phi, \ \ \ \ \  \varphi -\e  \psi = \delta  \phi.\]
Let $\Psi_{\e}=\Psi_{\e}(x, y)$ and $\bar{\Psi}_{\e}=\bar{\Psi}_{\e}(x, y)$ be defined in (\ref{p1}) and (\ref{p2}) respectively. 
Then 
\[ \Psi_{\e}(x, y) = \Theta(x, y), \ \ \ \ \  \bar{\Psi}_{\e}(x, y)
= \delta \Theta (\delta x, y ).\]

By Theorem \ref{thmK=-1}, we obtain the following
\begin{cor}\label{prop4.1}
Let $\phi(y)$ be a Minkowski norm on $\R^n$. 
Let $\Theta(x, y)$ denote the Funk metric of $\phi$. Then for any constant
$\delta$, the following function 
\be
F (x, y) := {1\over 2} \Big \{ \Theta(x, y) -\delta \Theta (\delta x, y )\Big \} \label{F1} 
\ee
is a projective Finsler metric on its domain with ${\bf K}=-1$ and its projective factor $P = {1\over 2} F^{-1} F_{x^k} y^k$ is given by
\be
P (x, y)
= {1 \over 2} \Big \{\Theta (x, y)
+\delta  \Theta (\delta  x, y )    \Big \}.\label{P1}
\ee
\end{cor}

\bigskip
\begin{ex}{\rm Take a look at  the special case when $\phi(y) = |y|$, 
the Funk metric  on $\Bbb B^n$ is given by
\be
 \Theta(x, y)= { \sqrt{|y|^2 - 
(|x|^2|y|^2 - \langle x, y \rangle^2) }+ \langle x, y \rangle \over 1- |x|^2}, \ \ \ \ \ y\in T_x \Bbb B^n.\label{Funkm}
\ee
Let
\begin{eqnarray}
F(x, y) & = & {1\over 2} \Big \{ { \sqrt{|y|^2 - 
(|x|^2|y|^2 - \langle x, y \rangle^2) }+ \langle x, y \rangle \over 1- |x|^2}\nonumber\\
&& - \delta { \sqrt{|y|^2 -\delta^2 
(|x|^2|y|^2 - \langle x, y \rangle^2) }+\delta  \langle x, y \rangle \over 1- \delta^2|x|^2}
\Big \}.
\end{eqnarray}
Clearly, $F$ is positively complete on
$\Bbb B^n(1)$  if $|\delta | < 1$ .
By Corollary \ref{prop4.1}, we know that
$F$ is projective with constant curvature ${\bf K}=-1$.
Note that when $\delta =-1$, $F$ is just the Klein metric on $\Bbb B^n(1)$.}
\end{ex}

\bigskip
Let $\phi(y)$ be an arbitrary Minkowski norm on $\R^n$ and $\Theta (x, y)$  denote the Funk metric of $\phi$. 
For   a constant vector 
$a\in \R^n$, let 
\[ \psi(y):= {1\over 2} \Big ( \phi(y)+ \langle a, y \rangle\Big ),
\ \ \ \ \ \varphi(y) := {1\over 2} \Big ( \phi(y)-\langle a, y\rangle \Big ),\]
such that 
\[ \varphi(y)+\psi(y) = \phi(y),
\ \ \ \ \ \varphi(y)-\psi(y) =- \langle a, y 
\rangle.\]
Let $\Psi_1=\Psi_{1}(x, y)$ and $\bar{\Psi}_1=\bar{\Psi}_{1}(x, y)$ be the function defined by (\ref{p1}) and (\ref{p2}) with $\e=1$, respectively.
We have
\[ \Psi_1(x, y) = \Theta(x, y), \ \ \ \ \ 
\bar{\Psi}_1 (x, y) = -{\langle a, y\rangle \over 1+\langle a, x\rangle }.
\]
By Theorem \ref{thmK=-1}, the following function 
\be
 F(x, y)= {1\over 2}\Big \{ \Theta(x, y) + {\langle a, y \rangle 
\over 1 + \langle a, x \rangle }\Big \} \label{bF}
\ee
is a projective Finsler metric  with ${\bf K}=-1$ and its projective factor is given by
\be
P(x, y) = {1\over 2}\Big \{ \Theta(x, y) - {\langle a, y \rangle 
\over 1 + \langle a, x \rangle }\Big \}.
\ee
We obtain the following 

\begin{cor}\label{prop4.2}  Let $\Theta(x, y)$ be the Funk metric on a strongly convex domain ${\cal U}\subset \R^n$ and $a$  an arbitrary constant vector. Let
\[ F(x, y):= {1\over 2} \Big \{ \Theta (x, y) +  {\langle a, y \rangle 
\over 1 + \langle a, x \rangle }\Big \}.
\]
 $F$  has the following properties
\ben
\item[(a)] ${\bf K}= - 1$;
\item[(b)] ${\bf E} = {1\over 2} (n+1) $;
\item[(c)] $F$ has straight geodesics. 
\een
\end{cor}
{\it Proof}: By the above argument, we know that $F$ is projective 
with ${\bf K}=-1$. Moreover the geodesic coefficients are in the form
$G^i = Py^i$. A direct computation gives 
\[ {\pa G^m\over \pa y^m} = (n+1) P .\]
Thus
\[ E_{ij} ={1\over 2} {\pa^2 \over \pa y^i \pa y^j} \Big [ {\pa G^m\over \pa y^m}  \Big ] =  {n+1\over 2} P_{y^iy^j}.\]
 By (\ref{Funk}) again, one immediately obtains 
\[ P(x, y)  = {1\over 2} \Big \{ \Theta (x, y)
- {\langle a, y \rangle 
\over 1 + \langle a, x \rangle } \Big \}.\]
Observe that
\begin{eqnarray*}
E_{ij} & = & {n+1 \over 4} \Big [ \Theta -{\langle a, y \rangle 
\over 1 + \langle a, x \rangle }\Big]_{y^iy^j} \\
& = & {n+1\over 4} \Big [ \Theta + {\langle a, y \rangle 
\over 1 + \langle a, x \rangle } \Big ]_{y^iy^j}\\
& = & {n+1\over 2} F_{y^iy^j}.
\end{eqnarray*}
This proves the corollary.
\qed

\bigskip
\begin{ex}{\rm Take a  look at the special case when $\phi(y)=|y|$. The Funk metric 
$\Theta(x, y)$ on the unit ball $\Bbb B^n(1)$ is given by (\ref{Funkm}).
Thus  for any constant vector $a \in \R^n$ with $|a| < 1$, the following function 
\be
F(x, y) ={1\over 2} \Big \{ { \sqrt{|y|^2 - 
(|x|^2|y|^2 - \langle x, y \rangle^2) }+ \langle x, y \rangle \over 1- |x|^2} +    {\langle a, y \rangle 
\over 1 + \langle a, x \rangle }\Big \}\label{eq50}
\ee
is a projective Finsler metric on $\Bbb B^n(1)$ with ${\bf K}=-1$ and ${\bf E}= {1\over 2}(n+1)$.
In \cite{Sh5}, we  prove that a Randers metric is projective with 
constant curvature if and only if it is locally Minkowskian or, up to a scaling,  isometric to the metric in (\ref{eq50}). }
\end{ex}

\section{${\bf K}=0$}

In this section, we are going to construct projective Finsler metrics $F(x, y)$ with ${\bf K}=0$, satisfying that 
 $F(0, y)=\psi(y)$ and $F_{x^k}(0, y)y^k =2 \psi(y) \varphi(y)$ for any given pair $\{ \psi(y), \varphi(y)\}$. We construct these examples using projective Finsler metrics $F_{\e}(x, y)$ with ${\bf K}=-\e^2$, satisfying $F_{\e}(0, y)=\psi(y)$ and $(F_{\e})_{x^k}(0, y)y^k = 2 \psi(y) \varphi(y)$.

\bigskip
Given  a constant $\e >0$, let  $\Psi_{\e}= \Psi_{\e}(x, y)$ and $\bar{\Psi}_{\e}=\bar{\Psi}_{\e}(x, y)$ denote the functions defined by 
(\ref{p1}) and (\ref{p2}) respectively. 
By Theorem \ref{thmK=-1}, the function
\be
F_{\e} (x, y):= {1\over 2\e} \Big \{  \Psi_{\e}(x, y) 
- \bar{\Psi}_{\e}(x, y) \Big \}
\ee
is a projective Finsler metric on its domain with ${\bf K}=-\e^2$ and 
its projector factor $P_{\e}:= {1\over 2}(F_{\e})^{-1} (F_{\e})_{x^k} y^k$ is given by
\be
P_{\e}(x, y) 
= {1\over 2} \Big \{  \Psi_{\e}(x, y) 
+ \bar{\Psi}_{\e}(x, y) \Big \}.
\ee
Let 
\begin{eqnarray*}
 \Psi(x, y) :& = & \lim_{\e\to 0^+} {\Psi_{\e} (x, y)-\Psi_0(x, y) \over \e},\\
\bar{\Psi}(x, y) :& = &  \lim_{\e\to 0^+} {\bar{\Psi}_{\e} (x, y)-\bar{\Psi}_0(x, y) \over \e}.
\end{eqnarray*}
Note that $\Psi_0(x, y)=\bar{\Psi}_0(x, y)$ and $\Psi_0=\Psi_0(x, y)$ is defined by
\be
\varphi\Big (y+ \Psi_0(x, y) \;x\Big ) = \Psi_0(x, y).\label{mmm}
\ee
Let
\be
F (x, y) := \lim_{\e\to 0^+} F_{\e}(x, y) = \Psi (x, y)\label{FK=0}
\ee
\be
P(x, y):= \lim_{\e\to 0^+} P_{\e}(x, y) = \Psi_0(x, y) \label{PK=0}
\ee
with $ F(0, y) = \Psi(0, y)= \psi(y)$ and $P(0, y)=\Psi_0(0, y)= \varphi(y)$. 

Differentiating (\ref{p1}) with respect to $\e$ at $\e=0$, we obtain
\begin{eqnarray}
\Psi(x, y) & = & {\psi \Big (y+\Psi_0(x, y) \; x\Big  ) \over 1- 
 \varphi_{y^k}\Big  (y+ \Psi_0(x, y) \; x\Big ) x^k }, \label{q1}\\
\bar{\Psi}(x, y) & =&  - {\psi \Big (y+\Psi_0(x, y) \; x\Big  ) \over 1- 
 \varphi_{y^k} \Big (y+ \Psi_0(x, y) \; x\Big ) x^k }. \label{q2}
\end{eqnarray}

We claim that $F(x, y):= \Psi(x, y)$ is a projective Finsler metric on its domain 
with ${\bf K}=0$ and its projective factor  is $P(x, y)= \Psi_0(x, y)$.

\begin{thm}\label{thmK=0}
Let $\psi(y)$ be a Minkowski norm on $\R^n$ and $ \varphi(y)$ a positively homogeneous function of degree one on $\R^n$. Let $\Psi_0=\Psi_0(x, y)$
de defined by (\ref{mmm}) and
$\Psi=\Psi(x, y)$  be  defined by (\ref{q1}). Then the function 
$F(x, y):=\Psi(x, y)$ is a projective Finsler metric on its domain with ${\bf K}=0$ and its projective factor $P(x, y)=\Psi_0(x, y)$.
Further, $F(0, y)=\psi(y)$ and $P(0, y)=\varphi(y)$.
\end{thm} 
{\it Proof}: Differentiating (\ref{mmm}) with respect to $x^k$ yields
\be
(\Psi_0)_{x^k}(x,y)
= {\varphi_{y^k}\Big  ( y + \Psi_0(x, y)x\Big ) \Psi_0(x, y)
\over 1- \varphi_{y^k}\Big  ( y + \Psi_0(x, y)x\Big )x^k }.\label{mmm1}
\ee
Contracting  (\ref{mmm1}) with $x^k$ yields
\[ \Psi_0(x, y) + (\Psi_0)_{x^k} (x, y) x^k 
= {\Psi_0(x, y) \over  1- \varphi_{y^k}\Big  ( y + \Psi_0(x, y)x\Big )x^k }.
\]
Thus we can express $\Psi(x, y)$  in (\ref{q1}) by 
\be
\Psi(x, y)={ \psi \Big (y+\Psi_0(x, y) \; x \Big )\over \Psi_0(x, y) } \Big \{
\Psi_0(x, y) + (\Psi_0)_{x^k}(x, y) x^k \Big \}.\label{q3}
\ee
Differentiating (\ref{mmm}) with respect to $y^k$ yields
\be
(\Psi_0)_{y^k} (x, y) = {\varphi_{y^k}\Big  ( y + \Psi_0(x, y)x\Big ) 
\over 1- \varphi_{y^k}\Big  ( y + \Psi_0(x, y)x\Big )x^k }.\label{mmm2}
\ee
It follows from (\ref{mmm1}) and (\ref{mmm2}) that 
\be
(\Psi_0)_{x^k}(x, y) = \Psi_0(x, y) (\Psi_0)_{y^k}(x, y).\label{PDEK=01}
\ee
Using  (\ref{q3}) and (\ref{PDEK=01})
one  can easily verify that 
\be
 \Psi_{x^k}(x, y) 
=(\Psi \Psi_0)_{y^k}(x, y) .\label{PDEK=0}
\ee
This is left to the reader. 
By Lemma \ref{lemkey}, $F(x, y)=\Psi(x, y)$ is a projective Finsler metric on its domain with ${\bf K}=0$ and its projective factor $P(x, y)=\Psi_0(x, y)$. \qed

\bigskip

\begin{cor}\label{prop5.1} Let $\phi(y)$ be a Minkowski norm on $\R^n$ and $\Theta(x, y)$ denote the Funk metric of $\phi(y)$. 
The following function 
\be
 F (x, y): =\Big \{ 1 +\langle a, x \rangle + {\langle a, y \rangle 
\over \Theta(x, y) }  \Big \}  \Big \{
\Theta(x, y) + \Theta_{x^k}(x, y) x^k \Big \}.
\ee
 is a projective Finsler metric on its domain with ${\bf K}=0$
and its projector factor $P= {1\over 2}F^{-1}F_{x^k}y^k$ is give by
\be
P(x, y) = \Theta(x, y).
\ee
\end{cor}
{\it Proof}: Take $\psi(y) =\phi(y)+ \langle a, y \rangle$ and $\varphi(y)= \phi(y)$. By the definition of $\Psi_0$, we have
\[ \Psi_0(x, y)=\Theta(x, y).\]
Observe that
\[ \phi ( y +\Psi_0(x, y) x)
=\phi ( y +\Theta(x, y) x) = \Theta(x, y).\]
Thus 
\[ \psi( y+ \Psi_0(x, y) x ) 
= \Theta(x, y) + \langle a, y\rangle + \langle a, x \rangle \Theta(x, y).\]
Therefore, the function $\Psi(x, y)$ in (\ref{q3}) is given by
\[ \Psi(x, y) = \Big \{ 1 +\langle a, x \rangle + {\langle a, y \rangle 
\over \Theta(x, y) }  \Big \}  \Big \{
\Theta(x, y) + \Theta_{x^k}(x, y) x^k \Big \}.\]
The corollary follows from Theorem \ref{thmK=0}. \qed

\bigskip

\begin{ex}{\rm
Take $\psi(y) = |y| + \langle a, y \rangle$ and $ \varphi(y) =|y|$. We obtain 
\begin{eqnarray}
F(x, y) & = & \Big \{ 1 +\langle a, x \rangle + { (1-|x|^2) \langle a, y \rangle \over \sqrt{|y|^2-(|x|^2|y|^2-\langle x, y \rangle^2) } +\langle x, y \rangle }\Big \} \nonumber\\
&& \times { \Big ( \sqrt{|y|^2-(|x|^2|y|^2-\langle x, y \rangle^2) } + \langle x, y \rangle \Big )^2\over 
(1-|x|^2)^2 \sqrt{ |y|^2-(|x|^2|y|^2-\langle x, y \rangle^2)} }. 
\end{eqnarray}
By Corollary \ref{prop5.1}, we know that $F(x, y)$ is a projective Finsler metric with ${\bf K}=0$ and its projective factor is given by
\be
P(x, y) = {\sqrt{|y|^2-(|x|^2|y|^2-\langle x, y \rangle^2) } +\langle x, y \rangle \over 1-|x|^2}.
\ee
}
\end{ex}

\bigskip

Let $\phi$ be an arbitrary  Minkowski norm on $\R^n$.
Take
\[ \psi(y):=   \phi(y) ,
\ \ \ \ \ \varphi(y):=  -\langle a, y \rangle.\]
We obtain 
\begin{eqnarray*}
 \Psi_0 (x, y) 
& = & - {\langle a, y\rangle \over 1+ \langle a, x\rangle },\\ 
\Psi(x, y) & = &  { \phi \Big ( (1+\langle a, x\rangle) y - \langle a, y\rangle x \Big ) \over \Big ( 1+ \langle a, x\rangle \Big )^2 } .
\end{eqnarray*}

By Theorem \ref{thmK=0}, we obtain the following
\begin{cor}\label{cor6.4} Let $\phi(y)$ be a Minkowski norm on $\R^n$ and $a\in \R^n$ a vector. The following function
\be
F(x, y) :=  { \phi \Big ( (1+\langle a, x\rangle) y - \langle a, y\rangle x \Big ) \over \Big ( 1+ \langle a, x\rangle \Big )^2 }\label{FFK=0}
\ee
is a projective Finsler metric on its domain with ${\bf K}=0$ and its projective factor $P(x, y)$ is given by
\be
P(x, y) = -{\langle a, y\rangle \over 1+ \langle a, x\rangle }.
\ee
\end{cor}

For the Finsler metric $F(x, y)$ in (\ref{FFK=0}), 
 the geodesic coefficients $G^i = P y^i$ are quadratic in $y\in \R^n$. Thus $F$ is a Berwald metric. It is well-known that any Berwald metric with ${\bf K}=0$ is locally Minkowskian. We conclude that
the Finsler metric in (\ref{FFK=0}) is locally Minkowskian.

\bigskip
\begin{ex}{\rm 
Take a Randers norm,  $\phi(y) = |y| + \langle b,  y \rangle$,
where $b\in \R^n$ is a vector with $|b| < 1$,  we obtain a projective Randers metric with ${\bf K}=0$,
\begin{eqnarray*}
F(x, y) & = &  { \sqrt{(1+\langle a, x \rangle )^2|y|^2
- 2 (1+\langle a, x\rangle )\langle a, y\rangle \langle x, y\rangle + \langle a, y\rangle^2 |x|^2  } \over ( 1+\langle a, x \rangle )^2 }\nonumber\\
&& + { (1+\langle a, x \rangle ) \langle b, y \rangle - \langle a, y \rangle \langle b, x \rangle \over ( 1+\langle a, x \rangle )^2 }.
\end{eqnarray*}
However, this metric must be locally Minkowskian. This fact also follows from  the main theorem in \cite{Sh5}. }
\end{ex}

\bigskip

Note that a projective Finsler metric $F(x, y):=\Psi(x, y)$ constructed in Theorem \ref{thmK=0} is reversible if and only if $F(0, y)=\psi(y)$ is reversible 
and $P(0, y)=\varphi(y)$ is anti-reversible. Thus there are lots of reversible projective Finsler metric with ${\bf K}=0$. 
We are going to show that any such metric is a Minkowski metric if it is complete. 
\begin{prop}\label{prop5.2}
Let $F$ be a complete reversible projective Finsler metric  on a strongly convex domain ${\cal U} \subset \R^n$. Suppose 
that $F$ has zero curvature ${\bf K}=0$. Then it is a Minkowski metric on 
$\R^n$.
\end{prop}
{\it Proof}: Let ${\bf v}\in T_x{\cal U}$ be an arbitrary vector
and $x (t):= x + t {\bf v}$. Let 
\[   F (t):= F \Big ( x(t), {dx\over dt}(t)\Big  ) ,
\ \ \ \ \ \ P(t):= P\Big  ( x(t), {dx\over dt}(t)\Big  ) .\]
By $ P = {1\over 2}F^{-1}F_{x^k}y^k$  and $ P^2 - P_{x^k}y^k =0$, we know that
$F(t)$ and $P(t)$ satisfy
\[ F'(t) = 2 F(t) P(t), \ \ \ \ \ P'(t) = P(t)^2.\]
We obtain 
\[ F(t) =  { F(0)\over \Big ( 1- P(0) t\Big  )^2 }, \ \ \ \ \ 
P(t) = { P(0) \over 1- P(0) t }.
\]
Assume that $P(0) >0$. 
Then 
\[\int_{-\infty}^0 F(t) dt = {F (0) \over P(0) } < \infty.\]
Thus $x(t)$ can not be extended to a $F$-geodesic defined on
$(-\infty, 0]$. If $P(0) < 0$, then 
\[ \int_0^{\infty} F(t) dt = - {F(0)\over P(0)} < \infty.\]
By assumption that $F$ is complete, we conclude that $P(0)= P(x,  {\bf v}) =0$ for any ${\bf v}\in T_x{\cal U}$. 
By (\ref{eq1}), we  obtain 
\[ F_{x^k} = (PF)_{y^k} = 0.\]
This implies that  $F(x, y) = F(0, y)$ is a Minkowski metric on $\R^n$.
\qed

\section{${\bf K}=1$}
In this section, 
we are going to construct projective Finsler metrics 
$F(x, y)$ with constant curvature ${\bf K}=1$ satisfying that 
$F(0, y)=\e \psi(y)$ and $F_{x^k}(0, y)y^k = 2 \psi(y)\varphi(y)$ for a  given pair $\{ \psi(y), \varphi(y)\}$ and $\e >0$.

Let $F(x,y)$ be a projective Finsler metric with ${\bf K}=1$ on an open neighborhood of the origin in $\R^n$. Suppose that $F(x, y)$ is $x$-analytic at $x=0$. 
By Theorem \ref{thm2}, $F_{\e}(x, y):=\e^{-1}F(x, y)$ 
and $P_{\e}(x, y)={1\over 2}(F_{\e})^{-1}(F_{\e})_{x^k}y^k$ are given by
\begin{eqnarray} 
F_{\e}(x, y) & = &  \sum_{m=0}^{\infty} {1\over m!}
{d^m\over dt^m} \Big [ \Psi_m (y + tx ) \Big ]|_{t=0},
\ \ \ \ \ \ y\in T_x\R^n,\label{eqKK=1}\\
P_{\e}(x, y) & = & \sum_{m=0}^{\infty} {1\over m!}
{d^m\over dt^m} \Big [ \Phi_m (y + tx ) \Big ]|_{t=0},\ \ \ \ \ \ y\in T_x\R^n,\label{eqPP=1}
\end{eqnarray}
where $\Psi_m(y)$ and $\Phi_m(y)$ are given by
\begin{eqnarray*}
\Psi_m  :& = &  { 1\over 2 (m+1)\e i} \Big \{ 
  \Big [ \varphi + i\e \psi  \Big ]^{m+1}
- \Big [  \varphi - i \e \psi \Big ]^{m+1}\Big \},\\
\Phi_m  : & = & {1\over 2 (m+1)}
\Big \{ \Big [\varphi + i\e \psi   \Big ]^{m+1} 
+ \Big [\varphi - i\e \psi   \Big ]^{m+1}\Big \},
\end{eqnarray*}
where $\psi(y)=F_{\e}(0, y)$ 
and $\varphi(y)=P_{\e}(0, y)$. 

\bigskip
Assume that $\psi(y)+i \e \psi(y)$ can be extended to  a complex-valued function on
$\Bbb C^n$. Define a complex-valued function  $H_{\e}= H_{\e}(x, y) =\Phi_{\e}(x, y)+i \Psi_{\e}(x, y)$ by
\be
H_{\e} = \varphi\Big ( y + H_{\e} x\Big  ) + i\e \psi\Big (y+ H_{\e} x\Big  ) .\label{DefH}
\ee
Differentiating (\ref{DefH}) with respect to $x^k$ and $y^k$ yields
\be
(H_{\e})_{x^k} = H_{\e} (H_{\e})_{y^k}. \label{HHH}
\ee
By (\ref{HHH}), we have
\be
(H_{\e})_{x^{i_1}\cdots x^{i_m}} (x, y)
= {1\over m+1} \Big [ (H_{\e})^{m+1}\Big ]_{y^{i_1}\cdots y^{i_m}} (x, y).\label{HHH1}
\ee
Note that $H_{\e}(0, y) = \varphi(y)+i\e \psi(y)$. Setting $x=0$ in (\ref{HHH1}) yields
\[(H_{\e})_{x^{i_1}\cdots x^{i_m}} (0, y)
= {1\over m+1} \Big [\Big  (\varphi +i\e \psi\Big  )^{m+1} \Big ]_{y^{i_1}\cdots y^{i_m}} (y).\]
Assume that $H_{\e}(x, y)$ is $x$-analytic. Then 
\be
H_{\e}(x, y)= 
\sum_{m=0}^{\infty} {1\over  (m+1)!} {d^m\over dt^m} \Big [ 
\Big (   \varphi(y+tx)+ i\e \psi(y+tx) \Big )^{m+1} \Big ]|_{t=0}.\label{PsiPsi}
\ee
One can directly verify that the function $H_{\e}(x, y)$ in (\ref{PsiPsi}) satisfies (\ref{HHH}). Thus we can also define $H_{\e}(x, y)$ by (\ref{PsiPsi}).

\bigskip

By (\ref{PsiPsi}), we can express $F_{\e}(x, y)$ in 
(\ref{eqKK=1})  and $P_{\e}(x, y)$ in (\ref{eqPP=1}) by 
\begin{eqnarray}
F(x, y) & = & \e F_{\e}(x, y)= {1\over 2i} \Big \{ H_{\e}(x, y)
- \overline{{H}_{\e} (x, y)} \Big \} = \Psi_{\e}(x, y),\label{eq52}\\
P(x, y) & = &  P_{\e}(x, y)= { 1 \over 2}
\Big \{ H_{\e} (x, y) 
+ \overline {H_{\e} (x, y) } \Big \}= \Phi_{\e}(x, y)\label{eq53}
\end{eqnarray}
with $F(0, y)=\e \psi(y)$ and $P(0, y)=\varphi(y)$. 

\begin{thm}\label{thmK=1}
Let $\psi(y)$ be a Minkowski norm on $\R^n$ and $\varphi(y)$ a positively homogeneous function of degree one on $\R^n$. Let $H_{\e} (x, y):=
\Phi_{\e}(x, y) + i \Psi_{\e} (x, y)$ denote the function defined in (\ref{DefH}) or (\ref{PsiPsi}).
Then 
$F(x, y):=  \Psi_{\e}(x, y)$ is a projective Finsler metric on its domain 
with ${\bf K}=1$ and its projective factor $P(x, y) = {1\over 2}F^{-1}F_{x^k}y^k$ is given by $P(x, y)= \Phi_{\e}(x, y)$. Further,
$F(0, y)=\e \psi(y)$ and $P(0, y)= \varphi(y)$. 
\end{thm}

\bigskip

Let
$\phi(y)$ be an arbitrary Minkowski norm on $\R^n$
and $\Theta(x, y)$ denote the Funk metric of $\phi$ defined in (\ref{Funkmetric}).
Let
\[ \psi(y):= \cos(\alpha) \phi(y),
\ \ \ \ \ \varphi(y):= \e\; \sin(\alpha) \phi(y),
\]
where $\alpha$ is an angle with $|\alpha | < \pi/2$.
Then
\[  \varphi(y)+i\e \psi(y) = i\e  e^{-i \alpha} \phi(y).\]
Assume that the  Funk metric 
$\Theta(x, y)$ can be extended to be an analytic function 
$\Theta(z, y)$ in $z\in {\cal U}\otimes \Bbb C \subset \Bbb C^n$. We have
\be
 H_{\e}(x, y) 
=i\e e^{-i\alpha} \Theta \Big (i \e e^{-i\alpha}x, y \Big ).\label{He}
\ee
By (\ref{eq52}) and (\ref{eq53}), we obtain the following
\begin{cor}{\rm Let $\phi(y)$ be a Minkowski norm on $\R^n$ and
 $\Theta(x, y)$ be the Funk metric of $\phi(y)$ on  ${\cal U}:=\{ y\in \R^n \ | \ \phi(y) <1\}$. Suppose that $\Theta(x, y)$ can be extended to a function 
$\Theta (z, y)$ on $({\cal U}\otimes \Bbb C) \times \R^n$. 
Define $H_{\e}(x, y)$ by (\ref{He}).
Then the following  function 
\be
F(x, y) := {\e \over 2} \Big \{ e^{-i\alpha} \Theta \Big ( i \e e^{-i\alpha}x, y \Big )+\overline{ e^{-i\alpha} \Theta \Big ( i\e  e^{-i\alpha}x, y \Big)} \Big \}
\ee
is a projective Finsler metric with ${\bf K}=1$ and its projector factor 
is given by 
\be
P(x, y) ={i\e \over 2} \Big \{ e^{-i\alpha} \Theta \Big ( i \e e^{-i\alpha}x, y \Big )- \overline{ e^{-i\alpha} \Theta \Big ( i\e  e^{-i\alpha}x, y \Big)}\Big \}.
\ee
Further, $F(0, y)= \e \cos(\alpha) \phi(y)$ and $P(0, y)=\e \sin(\alpha) \phi(y)$. }
\end{cor}

\bigskip
\begin{ex}{\rm 
Take a look at  the special case when
$\phi(y) = |y|$, i.e.,
$\psi(y) = \cos (\alpha) |y|$
and $\varphi(y) =\e \; \sin(\alpha) |y|$.
The Funk metric $\Theta(x, y)$ on the unit ball $\Bbb B^n$ is given by (\ref{Funkm}). 
We obtain
\begin{eqnarray*}
F(x, y) & = & {\e\over 2} \Big \{ e^{-2\alpha i}{\sqrt{ e^{2\alpha i} |y|^2 + \e^2 (|x|^2 |y|^2 -\langle x, y\rangle^2 )} + i \e \langle x, y \rangle \over 
1+ e^{-2\alpha i}\e^2 |x|^2 }\nonumber \\
&&  + e^{2\alpha i} {\sqrt{e^{-2\alpha i} |y|^2 +\e^2 (|x|^2 |y|^2 -\langle x, y\rangle^2 )} - i \e \langle x, y \rangle \over 
1 + e^{2\alpha i} \e^2 |x|^2 } \Big \},\\
P(x, y) & = & {i\e \over 2} \Big \{ e^{-2\alpha i}{\sqrt{ e^{2\alpha i} |y|^2 + \e^2 (|x|^2 |y|^2 -\langle x, y\rangle^2 )} + i \e \langle x, y \rangle \over 
1+ e^{-2\alpha i}\e^2 |x|^2 }\nonumber \\
&&  - e^{2\alpha i} {\sqrt{e^{-2\alpha i} |y|^2 +\e^2 (|x|^2 |y|^2 -\langle x, y\rangle^2 )} - i\e  \langle x, y \rangle \over 
1 + e^{2\alpha i}\e^2 |x|^2 } \Big \}.
\end{eqnarray*}

Both $F$ and $P$ take real values. Below we are going to get ride of $i$ in their expressions.
Let 
\begin{eqnarray*}
A : & = & \Big ( \cos(2\alpha) |y|^2 + \e^2 ( |x|^2 |y|^2 -\langle x, y\rangle^2 )\Big )^2 + \Big ( \sin(2\alpha) |y|^2 \Big )^2\\
B: & = & \cos(2\alpha) |y|^2 +\e^2 ( |x|^2 |y|^2 -\langle x, y\rangle^2 )\\
C: & = & \e \sin(2\alpha) \langle x, y \rangle\\
C': & = & \e \Big ( \cos(2\alpha) + \e^2 |x|^2 \Big ) \langle x, y\rangle ,\\
D: & = &\e^4 |x|^4 + 2 \cos(2\alpha)\e^2 |x|^2 + 1.
\end{eqnarray*}
Then for an angle $\alpha$ with $ 0\leq \alpha < \pi/2$,
\[
\sqrt{e^{2\alpha i} |y|^2 + \e^2 (|x|^2 |y|^2 -\langle x, y\rangle^2 )} = \sqrt{{\sqrt{A}+B\over 2}}+ i \sqrt{{\sqrt{A}-B\over 2}}.\]
\[
\sqrt{e^{-2\alpha i} |y|^2 + \e^2 (|x|^2 |y|^2 -\langle x, y\rangle^2 )} = \sqrt{{\sqrt{A}+B\over 2}}- i \sqrt{{\sqrt{A}-B\over 2}}.\]
We obtain 
\begin{eqnarray*}
F(x, y) & = & {\e\over D} \Big \{ \Big ( \cos(2\alpha) +\e^2|x|^2 \Big )
\sqrt{ \sqrt{A}+B \over 2} +\sin (2\alpha ) \sqrt{ \sqrt{A}-B \over 2} 
 +C \Big \}\\
P(x, y) & = & -{\e \over D}
\Big \{  \Big ( \cos(2\alpha) +\e^2|x|^2 \Big )\sqrt{ \sqrt{A}-B \over 2}-\sin (2\alpha) \sqrt{ \sqrt{A}+B \over 2}+C' \Big \}
.
\end{eqnarray*}

For an angle $\alpha$ with $ -\pi/2 < \alpha \leq 0$,
\begin{eqnarray*}
F(x, y) & = &  {\e \over D} \Big \{ \Big ( \cos(2\alpha) +\e^2 |x|^2 \Big )
\sqrt{ \sqrt{A}+B \over 2} - \sin (2\alpha ) \sqrt{ \sqrt{A}-B \over 2} +C\Big \}
,\\
P(x, y) & = &  -{\e \over D} \Big \{ \Big ( \cos(2\alpha) +\e^2 |x|^2 \Big )
\sqrt{ \sqrt{A}-B \over 2} + \sin (2\alpha ) \sqrt{ \sqrt{A}-B \over 2} +C'\Big \}.
\end{eqnarray*}
It is easy to verify that
\[
\Big ( \cos(2\alpha) +\e^2 |x|^2 \Big )
\sqrt{ \sqrt{A}+B \over 2} +|\sin (2\alpha )| \sqrt{ \sqrt{A}-B \over 2} 
=\sqrt{ { \sqrt{A}+B \over 2} D + C^2 },\]
\[
\Big ( \cos(2\alpha) +\e^2 |x|^2 \Big )
\sqrt{ \sqrt{A}-B \over 2} -|\sin (2\alpha )| \sqrt{ \sqrt{A}+B \over 2} 
=\sqrt{ { \sqrt{A}-B \over 2} D - C^2 }.\]
 $F$ and $P$ can be expressed by 
\begin{eqnarray}
F(x, y) & = & \e \Big \{ \sqrt{ { \sqrt{A}+B \over 2D}  + \Big ( {C\over D} \Big)^2 } + {C\over D}\Big \},\label{B55}\\
P(x, y) & =& -\e\Big \{  
\sqrt{ { \sqrt{A}-B \over 2D} - \Big ( {C\over D} \Big)^2 } + {C'\over D}\Big \}.
\label{B56}
\end{eqnarray}
One can verify that 
\[ \Xi:=P^2 - P_{x^k}y^k =F^2.\]
Thus $F$ has constant curvature ${\bf K}=1$.
In dimension two, one can verify that the Finsler metric 
$F$ in (\ref{B55}) is the  Bryant metric \cite{Br1}\cite{Br2}. 
}
\end{ex}

\noindent
Math Dept, IUPUI, 402 N. Blackford Street, Indianapolis, IN 46202-3216, USA.  \\
zshen@math.iupui.edu

\end{document}